\documentstyle[11pt,amssymb]{amsart}
\begin{document}

\title{Heat Kernel and Moduli Spaces II} 
\author{Kefeng Liu}

\date{November 1996 }
\maketitle

{\bf 1. Introduction.} 
In this paper we continue our study on the topology of the moduli spaces of flat bundles on a Riemann surface by using the heat kernels on compact Lie groups. As pointed out in [Liu], our method is very similar to the heat kernel proof of the Atiyah-Bott fixed point formula and the Atiyah-Singer index formula. In our case the local density is given by the Reidemeister torsion of the cochain complex of the simplest cell decomposition which gives the standard presentation of the fundamental group of the Riemann surface. The Reidemeister torsion is identified to the symplectic volume form, or the $L^2$ volume form, of the moduli space through a Poincare duality identity among the torsions. The global density is an infinite sum over all irreducible representations of the compact Lie group. This simple observation gives a complete answer to the questions about the intersection numbers of the cohomology classes of the moduli spaces of flat $G$-bundles for any compact semi-simple Lie!
!
!
 group $G$.  

The purposes of this note are several folds. The first is to extend the results in [Liu] to Riemann surfaces with several boundary components, to non-orientable Riemann surfaces and to the cases when the moduli spaces are singular; the second is to prove some general vanishing results about the characteristic numbers of the moduli spaces; the third is to apply our method to study the corresponding moduli spaces when $G$ is noncompact by using Hitchin's Higgs moduli spaces; the fourth is to introduce some invariants for knots and for $3$-manifolds by using the heat kernel method.  

Our approach to moduli spaces on Riemann surfaces is very flexible in applications. For example, we have only used the heat kernel of the Laplacian-Beltrami operator on the Lie group. One can instead use the heat kernel of more general biinvariant elliptic operators, or as in Section 7, consider to twist the integral of the pull-back of the heat kernel by the holonomy map by some general class functions from which we get the expressions for integrals of general functions on the moduli spaces in terms of heat kernel on the Lie group. This evaluates the integral of a general function on the moduli space with respect to the symplectic volume form, or the $L^2$ volume form. We hope to carry this out in a forthcoming paper.

Section $2$ is a review of the method of [Liu] in which Witten's formulas in [W1] are proved. The relationship between the Reidemeister torsion and the symplectic volume form is further clarified in this section. As corollaries we derive some general vanishing results about characteristic numbers of the moduli spaces. In Sect. 3 we derive integral formulas for moduli spaces on a Riemann surface with several boundary components; In Sect. 4, $L^2$ volume formulas for moduli spaces on non-orientable Riemann surfaces are obtained. In Sect. 5 we deal with singular moduli spaces. Our approach avoids some troubles appeared in [W1]. The volume formulas in these sections are due to Witten [W]. Note that our derivation of these volume formulas is much simpler than Witten's. In Sect. 6 we write down an integral formula on Hitchin's moduli spaces of Higgs pairs. This formula can be viewed as an integral formula for differential forms on the moduli spaces of flat bundles with non-compact s!
!
!
tructure groups. It may be usef

Note that the holonomy models for moduli spaces used in this paper are dual to the corresponding gauge theory models [AB], [Go], [Hu], [Hu1], as used in [AB], [W] and [W1]. Therefore it is the Reidemeister torsion, instead of the Ray-Singer torsion which is dual to the Reidemeister torsion, comes into the computations [RS]. Therefore in certain sense our approach is dual to the path-integral approach used by physicists [W1].

As remarked in [Liu], to get the Verlinde formula about the dimensions of nonabelian theta-functions on the moduli spaces from our approach, we only need to convert certain infinite sum into a finite sum. This should be a problem of Fourier transform and Poisson summation formula. But we are only able to do this by using the residue method of Szenes [Sz]. See the short note [LS] in which we have derived the general Verlinde formula from the formulas for the intersection numbers in this paper.

I would like to thank W. Zhang for explaining the Poincare duality about torsions to me. Many discussions with him and S. Chang are extremally helpful for the preparation of this note. I am also very grateful to the referee for his careful reading of the manuscript. This work has been partially supported by an NSF grant. 

\vspace{.1in}

\S 2. A review of the method.

\S 3. More boundary components.

\S 4. Non-orientable surfaces.

\S 5. Singular moduli spaces.

\S 6. Computation on Higgs moduli.

\S 7. Invariants for knots and $3$-manifolds.

\vspace{.2in}

{\bf 2.  A review of the method.}

We will use the same notations as in [Liu]. Let $G$ be a semi-simple simply connected compact Lie group and $T$ be a maximal torus. All of the discussions in this paper can be easily extended to the case when $G$ is not simply connected as in [Liu], Sect. 4. One simply divides our formulas by the order of the fundamental group of $G$, $\# \pi_1(G)$. Let ${\cal G},\ {\cal T}$ denote respectively the Lie algebras of $G$ and $T$. The basic idea in [Liu] is to use the holonomy model of the moduli space of flat principal $G$-bundles on an orientable Riemann surface $S$ of genus $g>1$ with boundary $\partial S$: 

$$ {\cal {M}}_{c}=f^{-1}  (c^{-1}) /Z_{c}$$ where $c\in T$ is a generic element, $Z_c\simeq T$ is the centralizer of $c$ and $f$ is the holonomy map:

$$f: \ G^{2g}\longrightarrow G$$
$$\ \ \ (x_1, y_1, \cdots x_g, y_g)\mapsto \prod_{j=1}^gx_j y_j x_j^{-1}y_j^{-1}.$$ The group $Z_c$ acts on $G^{2g}$ by conjugation 

$$\gamma(z)(x_1, \cdots, y_g)= (zx_1 z^{-1}, \cdots, z y_g z^{-1})$$ which induces an action on $f^{-1}(c^{-1})$. Equivalently we say ${\cal M}_c$ is the moduli space of flat connections on a principal $G$-bundle $P$, such that every element in ${\cal M}_c$ has fixed holonomy $c^{-1}$ around the boundary $\partial S$.

Let $$H(t, x, y)= \sum_{\lambda\in P_+}d_\lambda \chi_\lambda(xy^{-1})e^{-t p_c(\lambda)}$$ where $P_+\subset {\cal T}^*$ is the set of dominant integral weights of $G$ which is identified to the set of equivalent classes of irreducible representations of $G$, $p_c(\lambda)$ is the Casimir number, $d_\lambda$ is the dimension and $\chi_\lambda$ is the character of $\lambda$. It is known that $\frac{1}{\mbox{Vol}(G)} H(t, x, y)$ is the heat kernel of the Laplace-Beltrami operator on $G$ with respect to the bi-invariant metric $<\cdot, \cdot>$ induced from the Killing form [U], [Fe]. Let $\Delta^+$ denote the set, and $|\Delta^+|$ the number of positive roots. Recall that, in the induced metric on ${\cal T}^*$, $p_c(\lambda)=||\lambda+\rho||^2-||\rho||^2$ where 
$$\rho=\frac{1}{2}\sum_{\alpha\in \Delta^+} \alpha.$$  
Let $e\in G$ be the identity element and $\Gamma$ be the lattice in ${\cal T}$:
$$\Gamma=\{ H\in {\cal T};\ \mbox{exp}\, H=e\}.$$

Now consider the integral over $G^{2g}$ of the pull-back of $H(t, x, y)$ by $f$, 

$$I(t, c)= \int_{G^{2g}}H(t, f(h), c^{-1})\prod_{j=1}^g dx_j dy_j, \ \mbox{with}\ h=(x_1, y_1, \cdots, x_g, y_g)\in G^{2g}$$ where $dx_j, dy_j$ are the volume form on $G$ with respect to the bi-invariant metric. Let ${\mbox{Vol}}(G)$ and ${\mbox{Vol}}(T)$ denote the volumes of $G$ and respectively $T$ in the bi-invariant metric. By using the following well-known formulas for the characters of $G$ [BD]:

$$\int_G \chi_\lambda(wyzy^{-1}z^{-1})dz= \frac{\mbox{Vol}(G)}{d_\lambda}\chi_\lambda(wy)\chi_\lambda(y^{-1})$$ 
and 

$$\int_G \chi_\lambda(wy)\chi_\lambda(y^{-1})dy = \frac{\mbox{Vol}(G)}{d_\lambda}\chi_\lambda(w),$$we easily find 
 
$$I(t, c)= \mbox{Vol}(G)^{2g}\sum_{\lambda\in P_+}\frac{\chi_\lambda(c)}{d_\lambda^{2g-1}} e^{-tp_c(\lambda)}.$$

On the other hand, as $t\rightarrow 0$, $I(t, c)$ is localized to the neighborhood of $f^{-1}(c^{-1})$. By performing a simple Gaussian integral [Fo] and noting that the center $Z(G)$ of $G$ acts trivially on $f^{-1}(c^{-1})$, we get 

$$I(t, c)= \mbox{Vol}(G)\frac{\mbox{Vol}(T)}{\# Z(G)}\int_{{\cal {M}}_{c}}d\nu_c+O(e^{-\delta/t})$$ for any small postive number $\delta$ [Liu]. Here $\# Z(G)$ denotes the number of elements in $Z(G)$ and $d\nu_c$ is the torsion of the following chain complex [Fo], [Liu]:

$${\cal C}_c: \ \ \ 0\rightarrow {\cal Z}_{c}\stackrel{d\gamma}{\rightarrow}{\cal{G}}^{2g}\stackrel{df}{\rightarrow}{\cal{G}}\rightarrow 0 $$where ${\cal Z}_c$ denotes the Lie algebra of $Z_c$.

More precisely, ${\cal C}_c$ is the {\em cochain complex} $(C^*(S, \partial S, \mbox{ad}\, P), \partial)$ with coefficient in the adjoint bundle $\mbox{ad}\, P$ for the cell decomposition of $(S\, , \, \partial S)$ which gives the standard presentation of the fundamental group of the Riemann surface $S$ as shown in [Ma], pp38. Here $\partial$ denotes the coboundary operator. See the discussions in [Br], pp 59, [Ra], pp115, [Jo], [Go], [Gu], [Hu] and [Hu1] on such topic. One notes that $S$ is a $K(\pi, 1)$ manifold, therefore the cochain complexes for the group cohomology of $\pi_1(S)$ and the singular cohomology of $S$ coincide ([Br], pp59, [Hu], Sect. 4). Therefore we have 

$$d\nu_c= \tau(S, \partial S, \mbox{ad}\, P)$$ where $ \tau(S, \partial S, \mbox{ad}\, P)$ denotes the Reidemeister torsion of the cochain complex of the pair $(S, \partial S)$ with coefficient in $\mbox{ad}\, P$. 

Note that $H^0({\cal C}_c)$ and $H^2({\cal C}_c)$ both vanish for a generic element $c\in T$ and 

$$H^1({\cal C}_c)\simeq  H^1(S, \partial S, \mbox{ad}\, P)\simeq T{\cal M}_c.$$

Recall that Reidemeister torsion is equal to the Ray-Singer torsion by the theorem of Cheeger-Muller, which was originally conjectured by Ray-Singer. We refer the reader to [BZ], [Jo], [RS], [W], [Mi], [Fa] and [Lu] for discussions on torsions. Now we have the Poincare duality formula for torsions between the two comlexes, $(C^*(S, \partial S, \mbox{ad}\, P), \partial)$ and $(C^*(S, \mbox{ad}\, P), \partial)$,  which slightly generalizes the duality formula of [Mi] and was first shown to me by Zhang in this form: 
 
$$\tau(S, \partial S, \mbox{ad}\, P)\tau(S, \mbox{ad}\, P)=|| \mbox{det}\, H^1(S, \partial S,\mbox{ad}\, P)||^2_{L^2}$$where $\tau(S, \mbox{ad}\, P)$ denotes the torsion of $(C^*(S, \mbox{ad}\, P), \partial)$ and the right hand side is the induced $L^2$-metric on $ \mbox{det}\, H^1(S, \partial S,\mbox{ad}\, P)$. We also have the identity 

$$\tau(S, \mbox{ad}\, P)=\tau(S, \partial S, \mbox{ad}\, P)\tau(\partial S, \mbox{ad}\, P).$$
In our case, $\partial{S}$ is a circle and its torsion is easy to compute:
 $$\tau(\partial S, \mbox{ad}\, P)=|j(c)|^2$$ where 
$$|j(c)|^2=|\mbox{det}(1-\mbox{Ad}(c))|$$ is the square of the Weyl determinant. In fact, let $u\in Z(G)$ be an element in the center of $G$ and $c$ be an element near $u$, and $C\in {\cal T}$ be such that $c=u\, \mbox{exp}\, C$. Then    

$$j(c)=\prod_{\alpha\in \Delta^+}(e^{\sqrt{-1}\alpha(C)/2}-e^{\sqrt{-1}\alpha(C)/2}).$$ Here $\Delta^+$ is the set of positive roots of ${\cal G}$ with respect to ${\cal Z}_c\simeq {\cal T}$. Note that different choices of $C$ only change $j(c)$ by a plus or minors sign. We finally get 

$$|j(c)|d\nu_c=|j(c)|\tau(S, \partial S, \mbox{ad}\, P)=|| \mbox{det}\, H^1(S ,\partial S,\mbox{ad}\, P)||_{L^2}.$$

From the definition of the symplectic form $\omega_c$ on ${\cal{M}}_c$ ([AB], [W], [Ch1], [Je], or [Liu]), we easily see that the $L^2$-metric on $H^1(S,\partial S, \mbox{ad}\, P)\simeq T{\cal {M}}_{c}$,  which is the same as the Weil-Peterson metric, is the same as the metric induced by the symplectic form $\omega_{c}$ up to a normalization factor $4 \pi^2$. Therefore as norms on $\mbox{det}\, T{\cal M}_c$, we have the equality:

$$|j(c)|d\nu_c=(2\pi)^{2N_c}\frac{\omega_{c}^{N_c}}{N_c!}$$ where $N_c$ is the complex dimension of ${\cal M}_c$. Such equality was first used in [W], following [Jo], to study moduli spaces. By comparing the two expressions for $I(t, c)$, we have arrived at the formula:

$$ \ \ \int_{{\cal M}_c}e^{\omega_c}=
\# Z(G)\frac{\mbox{Vol}( G)^{2g-1}}{(2\pi)^{2N_c}\mbox{Vol}(T)}|j(c)|\sum_{\lambda\in P_+}\frac{\chi_\lambda(c)}{d_\lambda^{2g-1}}e^{-tp_c(\lambda)}+O(e^{-\delta/t}).$$

To derive the general formulas for the intersection numbers on ${\cal M}_u=f^{-1}(u^{-1})/G$ for $u^{-1} \in Z(G)$ a regular value of $f$, we use the relationship between the symplectic form $\omega_u$ on ${\cal {M}}_{u}$ and the symplectic form on ${\cal {M}}_{c}$ for $c$ a generic element near $u$. Then ${\cal{M}}_c$ is a fiber bundle over ${\cal{M}}_u$ with fiber $G/Z_c\simeq G/T$. Let 
$\pi: \ \ {\cal{M}}_c\rightarrow {\cal{M}}_u$ be the projection map, then we have 

$$\omega_c=\pi^* \omega_u +\nu_c$$
 where $\nu_c= <C, \Omega> +\theta_c$. Here recall that 
$C\in {\cal G}$ is given by $c=u\, \mbox{exp}\, C$, $2\pi\Omega$ is the curvature of the principal $G/Z(G)$-bundle $f^{-1}(u^{-1})\rightarrow {\cal{M}}_u$ and 
$\theta_c$, when restricted to each fiber of $\pi$, is the canonical symplectic form on $G/Z_c$ as given in [BGV], \S 7.5. See [Liu]. Note that ${\cal{M}}_c$ is actually euqivalent to $f^{-1}(u^{-1})/T$ [Liu]. One then has the simple push-forward formula:
$$\int_{{\cal{M}}_c}e^{\omega_c}= \int_{{\cal{M}}_u}e^{\omega_u}\pi_* e^{\nu_c}$$ and a family version of the Duistermaat-Heckman formula [KS], [Ch], [Liu]:

$$\pi_*e^{\nu_c}=\frac{\sum_{w\in W}\varepsilon(w)e^{<wC, \Omega>}}{\prod_{\alpha\in \Delta^+}(-\sqrt{-1}\alpha(\Omega))}$$where the sum is over the Weyl group $W$ and $\varepsilon(w)=\pm 1$ is the sign of $w\in W$. This formula can also be viewed as a family version of the localization formula in [BGV], Theorem 7.33. 

Choose an orthonormal basis of ${\cal T}$, $\{ H_1, \cdots, H_l\}$, then $C=x_1 H_1+\cdots +x_l H_l$. We use $(x_1, \cdots, x_l)$ as the coordinate of $C$ in ${\cal T}$. Let $p(x_1, \cdots ,x_l)$ be a Weyl-invariant polynomial. We then use the differential operator $p(\frac{\partial}{\partial x_1},\cdots, \frac{\partial}{\partial x_l})$ to act on both sides of the above formula which gives 

$$  \int_{{\cal{M}}_c} p(\sqrt{-1}\Omega)e^{\omega_c}  = \int_{{\cal{M}}_u} p(\sqrt{-1}\Omega)e^{\omega_u}\pi_* e^{\nu_c} \hspace{.6in}$$

$$\hspace{.4in}=\pm j(c)\# Z(G)\frac{\mbox{Vol}( G)^{2g-1}}{(2\pi)^{2N_c}\mbox{Vol}(T)}\sum_{\lambda\in P_+}\frac{\chi_\lambda(c)}{d_\lambda^{2g-1}}p(\lambda+\rho) e^{-tp_c(\lambda)}+O(e^{-\delta/t}).$$
Here we $p(\sqrt{-1}\Omega)$ is a differential form on ${\cal{M}}_u$, we also consider it as a form on ${\cal{M}}_c$ by pulling-back through the projection map $\pi$. One can also derive similar formulas even if $p$ is not Weyl-invariant.
 
By dividing both sides by $j(c)$ and taking limits, {\em first} $ t\rightarrow 0$, {\em then} $c\rightarrow u$, together with the formula [Liu]

$${\mbox{ lim}}_{c\rightarrow u}\frac{\pi_*e^{\nu_c}}{j(c)}=\pm\frac{{\mbox{Vol}}(G/T)}{(2\pi)^{{\mbox{dim}(G/T)}}},$$ we get the complete answers for the intersection numbers on ${\cal{M}}_u$, 

$$\int_{{\cal{M}}_u} p(\sqrt{-1}\Omega)e^{\omega_u} =\# Z(G)\frac{\mbox{Vol}( G)^{2g-2}}{(2\pi)^{2N_u}}\hspace{.7in}$$

$$\hspace{.6in} {\mbox{lim}}_{c\rightarrow u}{\mbox{lim}}_{t\rightarrow 0}\sum_{\lambda\in P_+}\frac{\chi_\lambda(c)}{d_\lambda^{2g-1}}p(\lambda+\rho)e^{-tp_c(\lambda)}\ \ \ (1).$$
More precisely we have obtained

(i). When $\mbox{deg}\, p\leq |\Delta^+|(2g-2)$, we have 

$$\int_{{\cal{M}}_u} p(\sqrt{-1}\Omega)e^{\omega_u} =\# Z(G)\frac{\mbox{Vol}( G)^{2g-2}}{(2\pi)^{2N_u}}\sum_{\lambda\in P_+}\frac{\Lambda_\lambda(u)}{d_\lambda^{2g-2}}p(\lambda+\rho)\ \ \ $$if the infinite sum on the right hand side converges. Here $\Lambda_\lambda(u)=\chi_\lambda(u)/d_\lambda$ is the induced character of the finite group $Z(G)$ in the representation $\lambda$ and $N_u$ is the complex dimension of ${\cal M}_u$. The $\pm$ sign is fixed by taking $p=1$.

(ii). When $\mbox{deg}\, p> |\Delta^+|(2g-2)$, 

$$\int_{{\cal{M}}_u} p(\sqrt{-1}\Omega)e^{\omega_u} =0$$ which follows from the Poisson summation formula, or Formula 4 in [Liu] applied inductively to each variable [W1].

From (ii) we easily obtain the following corollary:

\vspace{.1in}

{\bf Corollary 1:} {\em Let $P$ be a Pontryagin class of ${\cal{M}}_u$, if $\mbox{deg}\, P> |\Delta^+|(2g-2)$, then $\int_{{\cal M}_u}Pe^{\omega_u}=0$.}

\vspace{.1in}

Here $\mbox{deg}\, P$ means the degree of the Weyl-invariant polynomial correspponding to $P$. So as a differential form, the degree of $P$ is two times the degree used here. See the related work and conjectures in [Ne] for $G= SU(n)$ about the vanishing of Pontryagin classes. One can also use our method to derive vanishing results for Chern classes of ${\cal M}_u$. Note that here we have assumed that ${\cal{M}}_u$ is smooth. In the following, we will see that similar vanishing results also hold even when ${\cal{M}}_u$ is singular, or when $S$ has several boundary components.

Note that the order of taking limits is important in getting the general vanishing results in Corollaries 1, 2 and 3 in this paper. In [W], similar vanishing result on ${\cal{M}}_u$ is obtained under very restrictive condition on $u$. 

{\bf Remark 1.} To obtain the above corollary, we can go as follows instead which is much easier. Let 
$\pi(\lambda+\rho)=[\prod_{\alpha\in \Delta^+}<\alpha, \lambda+\rho>]$. It is easy to construct the differential operator $\pi(\partial)$ ([U], pp292) with respect to $C$ such that 

$$\pi(\partial)e^{(\lambda+\rho)(C)}=\pi(\lambda+\rho) e^{(\lambda+\rho)(C)}.$$
Let us denote 

$$V(c, t)=j(c)\sum_{\lambda\in P_+}\frac{\chi_\lambda(c)}{d_\lambda^{2g-1}} e^{-tp_c(\lambda)}.$$
Then we have, up to a constant $\pi(\rho)^{2g}$, 
$$\pi(\partial)^{2g}V(c, t)=j(c)\sum_{\lambda\in P_+}d_\lambda \chi_\lambda(c)e^{-tp_c(\lambda)}.$$
The proof of Formula 4 in [Liu] easily gives that, as $t$ goes to zero, the right hand side goes to zero exponentially for a generic $c$. This proves that, as $t$ goes to zero, the limit of $V(c, t)$ is a polynomial function in $C$ of degree at most $2g$. This method avoids the induction process [W1] and [Liu] for the vanishing results in Corollary 1 and in the following sections. Up to a constant, the operator $\pi(\partial)^2$ is the handle contraction operator in [Mo]. $\Box$

{\bf Remark 2.} Here are some corrections on [Liu]. In Formula 4 of [Liu], we need the polynomial $p(\lambda+\rho)=0$ if $\lambda$ is singular, that is if $\pi(\lambda+\rho)=0$. Equivalently we can require $p$ is of the form $q \, \pi(\lambda+\rho)$ where $q$ is any Weyl invariant or skew-invariant polynomial in $\lambda+\rho$. This condition was incidentally omitted in [Liu]. The following typos should be corrected too. In line 9, page 759 of [Liu], $C$ in the exponent should be $(C+H_u^0)$; in lines 6 and 9, page 758, the factor $2\pi$ in the denominators should be $(2\pi)^2$, since we used pairing $\alpha(C)$ which differs from [BGS] where inner product between $\alpha$ and the dual of $C$ is used for the volume formula, by a factor $2\pi$. Also in Lemma 4 and Lemma 4c, the power of $2\pi$ should be ${2N_u}$ and ${2N_c}$ respectively. $\Box$

\vspace{.2in}

{\bf 3. More boundary components.}

Now let $S$ be an orientable Riemann surface with $s$ disjoint boundary components. Let $c_1, \cdots, c_s\in G$ be $s$ generic elements. Here we assume $s\geq 2$, since $s=1$ case is already contained in Sect. 2. The standard presentation of the fundamental group of $S$ is given in [Ma], pp38, or any elementary text book on topology. We will build up the holonomy model for the moduli space of flat $G$-bundles in precisely the same way.

Let $f: \ G^{2g} \times G^{s-1} \rightarrow G$ be the holonomy map given by 

$$f(x_1, y_1,\cdots, x_g, y_g; z_2, \cdots z_s)= \prod_{j=1}^g x_j y_jx_j^{-1}y_j^{-1} \prod_{i=2}^s z_i c_i z_i^{-1}.$$
Let $Z_{c_j}$ be the centralizer of $c_j^{-1}$, then $Z_{c_1}\times Z_{c_2}\times \cdots \times Z_{c_s}$ acts on $G^{2g}\times G^{s-1}$ by 

$$\gamma (w_1, \cdots w_s)( x_1,\cdots, y_g; z_2, \cdots, z_s)= ( w_1x_1 w_1^{-1}, \cdots, w_1y_g w_1^{-1}; w_1z_2 w_2^{-1}, \cdots, w_1z_s w_s^{-1}).$$ 

The moduli space of flat connections on a flat $G$-bundle $P$ over $S$ with fixed holonomies $c_1^{-1}, \cdots, c_s^{-1}$ on the corresponding boundaries is just 

$${\cal{M}}_{c_1,\cdots, c_s}=f^{-1}(c_1^{-1})/Z_{c_1}\times \cdots \times Z_{c_s}.$$

Similar to [Liu] (see Sect. 2), we consider the integral 

$$I(t, c_1,\cdots c_s)= \int_{G^{2g}\times G^{s-1}}H(t, f(h), c^{-1}_1)\prod^g_{j=1} dx_j dy_j \prod^s_{i=2}dz_i.$$ 

One easily gets  

$$I(t, c_1,\cdots c_s)=\mbox{Vol}(G)^{2g+s-1}\sum_{\lambda\in P_+} \frac{\prod_{j=1}^s\chi_\lambda(c_j)}{d_\lambda^{2g-1+s}}e^{-tp_c(\lambda)}.$$

On the other hand, when $t\rightarrow 0$, we find 

$$I(t, c_1,\cdots, c_s)=\frac{\mbox{Vol}(G)}{\#Z(G)}\mbox{Vol}(T)^s\int_{{\cal{M}}_{c_1, \cdots,c_s}} d\nu_{c_1, \cdots, c_s}+O(e^{-\delta/t})$$ where $d\nu_{c_1, \cdots, c_s}$ is the torsion of 

$$0\rightarrow \oplus^s_{i=1} {\cal{Z}}_{c_i}\stackrel{d\gamma}{\rightarrow} {\cal{G}}^{2g}\oplus {\cal{G}}^{s-1}\stackrel{df}{\rightarrow} {\cal{G}}\rightarrow 0.$$
Again, this is precisely the cochain complex $(C^*(S, \partial S, \mbox{ad}\, P), \partial)$ with coefficient in $\mbox{ad}\, P$ of the cell decomposition of $(S, \partial S)$ which gives the standard presentation of the fundamental group of $S$ in [Ma], pp38.
 
By using the Poincare duality formula for the Reidemeister torsions in Section 2, we can identify $ d\nu_{c_1, \cdots, c_s}$ to the symplectic volume form on ${\cal{M}}_{c_1,\cdots, c_s}$:

$$ [\prod_{j=1}^s j(c_j)]\,d\nu_{c_1, \cdots, c_s}= (2\pi)^{2N}\frac{\omega^N_{c_1, \cdots, c_s}}{N!}$$ where $N$ is the complex dimension of ${\cal{M}}_{c_1,\cdots, c_s}$, and $\omega_{c_1, \cdots, c_s}$ is the symplectic form on $ {\cal{M}}_{c_1,\cdots, c_s}$. Here just note that the term $\prod_{j=1}^s j(c_j)^2$ is just the Reidemeister torion of the boundary $\partial S$ which is $s$ copies of $S^1$, and the above formula follows from the argument in Section 2.

Therefore we have arrived at the following volume formula for the moduli space of flat bundles on $S$:

$$\int_{{\cal{M}}_{c_1, \cdots, c_s}}e^{\omega_{c_1, \cdots, c_s}}=\frac{\#Z(G){\mbox{Vol}}(G)^{2g-2+s}}{(2\pi)^{2N}{\mbox{Vol}}(T)^s} [\prod_{j=1}^s j(c_j)] \sum_{\lambda\in P_+}\frac{\prod_{j=1}^s\chi_\lambda(c_j)}{d_\lambda^{2g-2+s}}e^{-tp_c(\lambda)}+O(e^{-\delta/t}).$$

As in Section 2, by introducing an extra puncture and taking Weyl-invariant derivatives we get the formulas for general intersection numbers on ${\cal{M}}_{c_1,\cdots, c_s}$. 

$$\int_{{\cal{M}}_{c_1, \cdots, c_s}} p(\sqrt{-1}\Omega)e^{\omega_{c_1, \cdots, c_s}}=\frac{\#Z(G){\mbox{Vol}}(G)^{2g-2+s}}{(2\pi)^{2N}{\mbox{Vol}}(T)^s}\hspace{.3in}$$

$$\hspace{.7in} [ \prod_{j=1}^s j(c_j)] \ {\mbox{lim}}_{t\rightarrow 0}\sum_{\lambda\in P_+}\frac{\prod_{j=1}^s\chi_\lambda(c_j)}{d_\lambda^{2g-2+s}}p(\lambda+\rho)e^{-tp_c(\lambda)}\hspace{.2in} (2)$$
where $p$ is a Weyl-invariant polynomial and $2\pi \Omega$ is the curvature of the universal principal $G$-bundle on ${\cal{M}}_{c_1, \cdots, c_s}$. This universal $G$-bundle is given by the restriction of the universal $G$-bundle on $S\times  {\cal{M}}_{c_1, \cdots, c_s}$ to $\{o\}\times {\cal{M}}_{c_1, \cdots, c_s}$ for some point $o\in S$. It is also the pull-back of the principal $G$-bundle $f^{-1}(c_1)\rightarrow {\cal{M}}_{c_2, \cdots, c_s}$ to ${\cal{M}}_{c_1, \cdots, c_s}$ by the projection map ${\cal{M}}_{c_1, \cdots, c_s}\rightarrow {\cal{M}}_{c_2, \cdots, c_s}$ of specializing $c_1$ to $e$.  

In particular for $\mbox{deg}\, p\leq (2g-2+s)|\Delta^+|$, we have  

$$\int_{{\cal{M}}_{c_1, \cdots, c_s}} p(\sqrt{-1}\Omega)e^{\omega_{c_1, \cdots, c_s}}=\hspace{1.1in}$$

$$\hspace{.6in}\frac{\#Z(G){\mbox{Vol}}(G)^{2g-2+s}}{(2\pi)^{2N}{\mbox{Vol}}(T)^s}[ \prod_{j=1}^s j(c_j)] \sum_{\lambda\in P_+}\frac{\prod_{j=1}^s\chi_\lambda(c_j)}{d_\lambda^{2g-2+s}}p(\lambda+\rho)\hspace{.2in}$$if the right hand side converges.

It is also straightforward to derive the following vanishing result:

\vspace{.1in}

{\bf Corollary 2:}  {\em Let $P$ be a Pontryagin class of $ {\cal{M}}_{c_1,\cdots, c_s}$, if $\mbox{deg}\, P> (2g-2+s)|\Delta^+|$, then $\int_{{\cal M}_{c_1, \cdots, c_s}}Pe^{\omega_{c_1,\cdots, c_s}}=0$.}

\vspace{.1in}

Note that here we have assumed $c_1, \cdots, c_s$ are generic elements in $G$, so that the moduli spaces in the above discussion are all smooth. But all the results can be extended to the singular case by following the discussions in Sect. 5.

{\bf Remark 3.} It is clear that to each puncture there corresponds a principal $T$-bundle. When we take derivatives with respect to the $c_j$'s, the curvatures of these $T$-bundles will appear in the integrand of Formula (2). In this way we get general formulas for the integrals of other generators of the cohomology groups of ${\cal{M}}_{c_1, \cdots, c_s}$, see [BR].$\Box$   

In the above discussions we can use instead the holonomy model of the moduli space corresponding to different cell decomposition of $S$ as given in [Ma], pp38. Then the corresponding cochain complex is

$$0\rightarrow {\cal{G}}\oplus\oplus _{i=1}^s {\cal{Z}}_{c_i}\stackrel{d\gamma}{\rightarrow} {\cal{G}}^{2g}\oplus {\cal{G}}^s\stackrel{df}{\rightarrow} {\cal{G}}\rightarrow 0.$$
In this case the holonomy map $f$ is given by 

$$f:\ G^{2g}\times G^s\rightarrow G$$
$$\ \ \ (x_1, y_1,\cdots, x_g, y_g; z_1, \cdots, z_s)\mapsto \prod_{j=1}^g x_jy_j x_j^{-1}y_j^{-1}\prod_{i=1}^s z_ic_iz_i^{-1}$$
and the action $\gamma$ of $G\times Z_{c_1}\times \cdots \times Z_{c_s}$ on $G^{2g}\times G^s\rightarrow G$ is  

$$\gamma(w; w_1,\cdots, w_s)(x_1,\cdots, y_g; z_1, \cdots, z_s)=
(wx_1 w^{-1}, \cdots, wy_gw^{-1}; wz_1w_1^{-1}, \cdots, wz_sw_s^{-1}) .$$

As we can see from the above discussions that the presentations of the fundamental group of the Riemann surfaces completely determine our ways of computations on the moduli spaces of flat $G$-bundles. But all the different ways give the same answer.

\vspace{.2in}

{\bf 4. Non-orientable surfaces.}

Depending on odd or even copies of $RP^2$ glued together, the presentation of the fundamental group of a closed non-orientable Riemann surface $S$ is given as follows:

\vspace{.1in} 

(i). $\pi_1(S)$ has $2k+1$ generators $\{ x_j,\ y_j, \varepsilon; \ j=1, \cdots, k\} $ with one relation:

$$\prod_{j=1}^k [x_j y_jx_j^{-1}y_j^{-1}] \varepsilon^2;$$

(ii). $\pi_1(S)$ has $2k+2$ generators $\{ x_j, \ y_j, \ z,\ \varepsilon; \ j=1, \cdots ,k\}$, with a single relation:

$$ \prod_{j=1}^k [x_j y_jx_j^{-1}y_j^{-1}] z\varepsilon z^{-1}\varepsilon.$$
See [Ma], pp135. We could also use the presentations in [Ma], pp38 which work in the same way. Here we choose this one to keep the computations and notations compatible with Section 2. The above presentations clearly indicate the way to compute on the corresponding moduli spaces.

Let $u\in Z(G)$ be an element in the center of $G$. We first consider Case (i) in which the holonomy map is 

$$f:\ \ G^{2k}\times G\rightarrow G$$

$$\ \ (x_1, \cdots, y_k; \varepsilon)\mapsto \prod_{j=1}^k x_j y_jx_j^{-1}y_j^{-1} \varepsilon^2.$$    
The conjugate action $\gamma$ of $G$ on $G^{2k}\times G$ is just

$$\gamma(w)(x_1, \cdots, y_k; \varepsilon)=(wx_1w^{-1}, \cdots, wy_kw^{-1}; w\varepsilon w^{-1}).$$

Assume $u^{-1}$ is a regular value of $f$, then the moduli space is easily obtained as  

$${\cal{M}}_u =f^{-1}(u^{-1})/G.$$
The real dimension of ${\cal M}_u$ is $(2k-1)\mbox{dim}\, G$. Consider the integral :

$$I(t, u)=\int_{G^{2k}\times G} H(t, f(h), u^{-1})\prod^k_{j=1}dx_j dy_j d\varepsilon$$where $h=(x_1, \cdots, y_k; \varepsilon)\in G^{2k}\times G$. Similar to the computations in Sect. 2, together with the formula 

$$\int_G\chi_\lambda(\varepsilon^2)d\varepsilon=f_\lambda \cdot \mbox{Vol}(G)$$ where $f_\lambda=1$ if $\lambda$ has a real structure, $-1$ if $\lambda$ has a quarternionic structure and $0$ otherwise. We get 

$$I(t, u)={\mbox{Vol}}(G)^{2k+1}\sum_{\lambda\in P_+}f_\lambda \frac{\chi_\lambda(u)}{d_\lambda^{2k-1}}e^{-t p_c(\lambda)}.$$When $t\rightarrow 0$, we have 

$$I(t, u)=\frac{{\mbox{Vol}}(G)^2}{\#Z(G)}\int_{{\cal{M}}_u}d\nu_u+ O(e^{-\delta/t})$$ for any small positive number $\delta$. Here $d\nu_u$ is the torsion of 

$$0\rightarrow {\cal G}\stackrel{d\gamma}{\rightarrow} {\cal{G}}^{2k}\oplus {\cal G}\stackrel{df}{\rightarrow} {\cal{G}}\rightarrow 0$$
which is the cochain complex $(C^*(S, \mbox{ad}\, P), \partial )$ associated to the cell decomposition of the presentation of the fundamental group of $S$ as given in [Ma], pp135. There is no symplectic structure on ${\cal M}_u$, but the form 

$$\frac{1}{4\pi^2}\int_{S}<a,\ b>,\ {\mbox{for}}\ a,\ b \in H^1(S, {\mbox{ad}}\, P)$$ which induces the symplectic form when $S$ is orientable, still induces the $L^2$ volume form on ${\cal M}_u$. The Poincare duality argument in Sect. 2 shows that $d\nu_u$ is identified to this $L^2$ volume form, or the Weil-Peterson volume form $dv_{L^2}$ up to a normalization by $4\pi^2$:

$$d\nu_u= (2\pi)^{N} dv_{L^2}$$ where $N=(2k-1)\mbox{dim}\,G$ is the real dimension of ${\cal{M}}_u$. Let $ V_{L^2}({\cal M}_u)$ denote the $L^2$ volume of ${\cal M}_u$. Let $t\rightarrow 0$, we easily get Witten's formula [W] from the above two expressions for $I(t, u)$:

$$ V_{L^2}({\cal M}_u)=\# Z(G)\frac{{\mbox{Vol}}(G)^{2k-1}}{(2\pi)^{N}} \sum_{\lambda\in P_+}f_\lambda \frac{\Lambda_\lambda(u)}{d_\lambda^{2k-2}}.$$

For Case (ii) of the fundamental group, one only needs to consider the holonomy map:

$$f:\ G^{2k}\times G\times G\rightarrow  G$$

$$\ \ \  (x_1, \cdots, y_k; z; \varepsilon)\mapsto \prod_{j=1}^k x_jy_j x_j^{-1}y^{-1}_j z \varepsilon z^{-1}\varepsilon$$ with the action $\gamma$ of $G$ on each factor by conjugation. The corresponding moduli space ${\cal{M}}_u$ for a regular value $u^{-1}\in Z(G)$ is simply given by 

$${\cal{M}}_u=f^{-1}(u^{-1})/G.$$ Then the computations are completely the same as in Case (1). We get the $L^2$ volume formula 

$$V_{L^2}({\cal M}_u)=\# Z(G)\frac{{\mbox{Vol}}(G)^{2k}}{(2\pi)^{\mbox{dim}\, {\cal M}_u}} \sum_{\lambda\in P_+}f_\lambda \frac{\Lambda_\lambda(u)}{d_\lambda^{2k}}$$where $\mbox{dim}\, {\cal M}_u=2k\mbox{dim}\, G$.

We can also combine with the computations in Sect. 3 to get volume formulas of the moduli spaces of flat bundles over a non-orientable Riemann surface with several boundary components. The interested reader may enjoy the computation by working out the detail.

\vspace{.2in} 

{\bf 5. Singular moduli spaces.} 
                                                       
Let $S$ be an orientable Riemann surface of genus $g>1$. In [Go1] and [Hu] it is shown that, when ${\cal{M}}_u$ with $u\in Z(G)$ is singular, it is a stratified symplectic manifold. Each stratum is a smooth symplectic manifold with induced symplectic form from the standard symplectic form on $H^1(S, \mbox{ad}\, P)$. It has a Zariski open dense stratum ${\cal M}_u^o$. The symplectic volume on each stratum is finite [Hu]. The codimension of the singularity is at least $1$, actually larger than $4$ except when $g=2$ and $G=SU(2)$ [F]. Let $c\in G$ be a generic point near $u$, such that the moduli space ${\cal M}_c$ is smooth. We assume the complex dimension of ${\cal{M}}_u$ is still $N_u$ . 

The basic idea in this section is to use ${\cal{M}}_c$ as a smooth model for ${\cal{M}}_u$. We lift the computations on ${\cal{M}}_u$ to ${\cal{M}}_c$. Technically nothing is substantially different from the smooth case.

For a Weyl-invariant polynomial $p$, from our previous formulas we have 

$$\frac{1}{j(c)}\int_{{\cal{M}}_{c}}p(\sqrt{-1}\Omega)e^{\omega_c}=\frac{\#Z(G){\mbox{Vol}}(G)^{2g-1}}{(2\pi)^{2N_c}{\mbox{Vol}}(T)} \sum_{\lambda\in P_+}\frac{\chi_\lambda(c)}{d_\lambda^{2g-1}}p(\lambda+\rho)e^{-tp_c(\lambda)}+\frac{1}{j(c)}O(e^{-\delta/t})$$where $2\pi\Omega$ is the curvature of the universal $G$-bundle on ${\cal{M}}_{c}$. See [BR] for some basic facts about universal bundles on ${\cal{M}}_c$. Motivated by the smooth case, we consider the limit

$$\mbox{lim}_{c\rightarrow u}\frac{1}{j(c)}\int_{{\cal{M}}_{c}}p(\sqrt{-1}\Omega)e^{\omega_c}.$$
We expect such limits to give us all needed information about integrals on ${\cal{M}}^o_u$. To justify this let us understand how much information about ${\cal{M}}_u$ is encoded in this limit.  

Let $P^o\rightarrow {\cal{M}}_u^o$ be the universal $G$-bundle on ${\cal{M}}_u^o$ and $2\pi\Omega^o$ be its curvature. Recall that $G$ acts on $f^{-1}(u^{-1})$ by conjugation, $P^o$ can be identified as the union of the principal $G$-orbits in $f^{-1}(u^{-1})$ ([Hu], Theorem 5). It is clear that $P^o/G={\cal{M}}_u^o$ and the quotient $P^o/T$ is equivalent to a Zariski open dense stratum ${\cal{M}}_c^o$ in $ {\cal{M}}_c$. The projection $\pi: \  {\cal{M}}_c^o\rightarrow {\cal{M}}_u^o$ is a fibration with fiber $G/T$. In fact $\pi$ can be extended to a dominant map from ${\cal{M}}_c$ to ${\cal{M}}_u$, but this is not needed in the following. The codimension of ${\cal{M}}_c-{\cal{M}}_c^o$ is also at least $1$,

The pull-back $\pi^*p(\sqrt{-1}\Omega^o)$ to ${\cal{M}}_c^o$ has an extension to ${\cal{M}}_c$ which is given by $p(\sqrt{-1}\Omega)$ by the functoriality of these classes [AB], [BR]. In fact both of them are equal to the pull-back of the same universal class on the classifying space of $G$ by the natural classifying maps. See [BR], Prop. 3.1. Equivalently one can easily see that the restriction of the universal $G$-bundle to ${\cal{M}}_c^o$ is the same as $\pi^* P^o$.

Since the codimension of ${\cal{M}}_c-{\cal{M}}_c^o$ is at least $1$, we have 

$$\int_{{\cal{M}}_c}p(\sqrt{-1}\Omega)e^{\omega_c}=\int_{{\cal{M}}_c^o}\pi^*p(\sqrt{-1}\Omega^o)e^{\omega_c}=\int_{{\cal{M}}_u^o}p(\sqrt{-1}\Omega^o)e^{\omega_u}\pi_*e^{\nu_c}.$$ 
By applying our previous argument to the $G/T$-fibration $\pi:\ {\cal{M}}_c^o\rightarrow{\cal{M}}_u^o$, we obtain

$$\frac{\mbox{Vol}(G/T)}{(2\pi)^{\mbox{dim}( G/T)}}\int_{{\cal{M}}_u^o}p(\sqrt{-1}\Omega^o)e^{\omega_u}=\mbox{lim}_{c\rightarrow u}\frac{1}{j(c)}\int_{{\cal{M}}_c}p(\sqrt{-1}\Omega)e^{\omega_c}$$

$$\hspace{.6in}= \# Z(G)\frac{{\mbox{Vol}}(G)^{2g-1}}{(2\pi)^{2N_c}\mbox{Vol}(T)}{\mbox{lim}}_{c\rightarrow u}{\mbox{lim}}_{t\rightarrow 0}\sum_{\lambda\in P_+}\frac{\chi_\lambda(c)}{d_\lambda^{2g-1}}p(\lambda+\rho) e^{-tp_c(\lambda)}\hspace{.2in} (3).$$
Here we have used our previous formula about symplectic forms in Section 2 and [Liu]:
$$\omega_c= \pi^*\omega_u +\nu_c$$ togehter with the identity 

$$\mbox{lim}_{c\rightarrow u}\frac{\pi_* e^{\nu_c}}{j(c)}=\pm \frac{\mbox{Vol}(G/T)}{(2\pi)^{\mbox{dim}( G/T)}}.$$ 
As pointed out in [Liu], the above formulas are simple corollaries of the local model theorem for symplectic manifold which still holds when the base symplectic manifold is open [GS]. In fact $2\pi\nu_c=< C, d\theta>$ where $\theta$ is the connection $1$-form on $P^o$ whose curvature is $2\pi \Omega^o$ and $C$ is the Lie algebra element of $c$ with $c=u\,{\mbox{exp}}\,C$. By using the curvature formula 

$$2\pi \Omega^o=d\theta+\frac{1}{2}[\theta, \theta],$$ we can write 

$$2\pi\nu_c=<C, 2\pi \Omega^o>-\frac{1}{2}<C, [\theta, \theta]>.$$ Up to a factor $2 \pi$, the last term of this identity is just $\theta_c$ of Section 2 which induces the canonical symplectic form on each fiber of $\pi:\ {\cal{M}}^o_c\rightarrow {\cal{M}}_u^o$. See [Liu], pp757 and Section 2. 

Formula (3) contains two parts: 

(i). The integral $\int_{{\cal{M}}_u^o}p(\sqrt{-1}\Omega^o)e^{\omega_u}$ is well-defined. By considering ${\cal{M}}_c$ as ``the desingularization'' of ${\cal{M}}_u$, we may view the first identity of (3) as the definition of this integral.

(ii). We have precise formulas for the integals on ${\cal{M}}_u^o$ given by the limits of the infinte sum in the second equality of (3). 

One should note that the order of taking limits in (3) is important, {\em first} $t$, {\em then} $c$. The reason is that, if we take the $c$ limit first, the term $O(e^{-\delta/t})/j(c)$ may go to infinity, and this may also cause the divergence of the infinite sum in (3). 

In particular we have 

$$\int_{{\cal{M}}_u^o}p(\sqrt{-1}\Omega^o)e^{\omega_u}=\# Z(G)\frac{\mbox{Vol}( G)^{2g-2}}{(2\pi)^{2N_u}}\sum_{\lambda\in P_+}\frac{\Lambda_\lambda(u)}{d_\lambda^{2g-2}}p(\lambda+\rho)$$
if the infinite sum on the right hand side converges.

 We also record the following obvious corollary which is considered as the corresponding vanishing result for the singular moduli spaces:

\vspace{.1in}

{\bf Corollary 3.}  {\em Let $p$ be a Weyl-invariant polynomial. If $\mbox{deg}\, p>(2g-2)|\Delta^+|$, then }

$$\int_{{\cal{M}}^o_u}p(\sqrt{-1}\Omega^o)e^{\omega_u}=\mbox{lim}_{c\rightarrow u}\frac{1}{j(c)}\int_{{\cal{M}}_{c}}p(\sqrt{-1}\Omega)e^{\omega_c}=0.$$

\vspace{.1in}

The results in this section were speculated in [W1]. In dealing with singular moduli spaces, some trouble appears with the fractional power of the parameter $\varepsilon$ in [W1]. This is avoided by our approach. One can easily extends the results to the case when the Riemann surface $S$ has several boundary components. Also the method can be used to extend the volume formulas in Sect. 4 to the singular case.

{\bf Remark 4.} In many cases, we can work directly over ${\cal{M}}_c$ and avoid the process of taking the $c$-limit. In [LS], the idea of using ${\cal{M}}_c$ as a smooth model for ${\cal{M}}_u$ was pursued to compute the Riemann-Roch number of the line bundle $L_u$ on ${\cal{M}}_u$ with $c_1(L_u)=\omega_u$. In fact, since the dominant map $\pi:\ {\cal{M}}_c\rightarrow  {\cal{M}}_u $ of specializing $c$ to $u$ has generic fiber $G/T$ and the singularity has codimension at least $1$, we have 

$$\chi({\cal{M}}_u, L_u^k)=\chi({\cal{M}}_c, \pi^*L_u^k)$$ for any positive integer $k$. Here $\chi$ denotes the Riemann-Roch number which is the alternating sum of the dimensions of the cohomology groups of the corresponding line bundle. Note that ${\cal{M}}_c$ is smooth and $\pi^*L_u$ extends to a line bundle on ${\cal{M}}_c$. By Atiyah-Singer index theorem we have 

$$\chi({\cal{M}}_c, \pi^*L_u^k)=\int_{{\cal{M}}_c}\mbox{Td}({\cal{M}}_c)\mbox{ch}\,\pi^*L_u^k=\int_{{\cal{M}}_c}\mbox{Td}({\cal{M}}_c)e^{k\pi^*\omega_u}$$
and in principle the right hand side can be computed by our formulas in Section 3, once one knows the expression of the Todd class $\mbox{Td}({\cal{M}}_c)$ in terms of the generators of the cohomology groups of ${\cal{M}}_c$. This is currently being worked out in [LS]. $\Box$

\vspace{.2in}

{\bf 6. Computation on Higgs moduli.}

Hitchin's Higgs moduli space in [Hi] can be considered as the moduli space of flat bundles on an orientable compact Riemann surface $S$ with reductive structure groups. Therefore we have a good model to extend our results to noncompact Lie groups.

Let $G$ be a compact Lie group and $P$ be a principal $G$-bundle over the Riemann surface $S$. A Higgs pair is a connection $A$ on $P$ and a section $\Phi$ of $K\bigotimes \mbox{ad}\, P_ C$ where $K$ is the cotangent bundle of $S$ and ${\mbox{ ad}}\, P_C$ is the complexified adjoint bundle, such that 

$$F_A= -[\Phi, \, \Phi^*],\ \ \ \bar{\partial}_A \Phi=0$$ where $\bar{\partial}_A$ is the $(0, 1)$-part of the covariant derivative of the connection $A$. 
 
Note that a Higgs pair gives us a flat connection on $\mbox{ad}\, P_C $, therefore a representation of $\pi_1(S)$ in the reductive Lie group $G_C $ [Hi]. 

Assume the genus of $S$ is larger than 1. The moduli space of Higgs pairs ${\cal{M}}_H$ is, roughly speaking, like the cotangent bundle of ${\cal{M}}_u$. For simplicity, we will only deal with the situation that ${\cal{M}}_H$ is smooth. Then ${\cal{M}}_H$ is a smooth symplectic manifold and there is a circle action on ${\cal{M}}_H$ with a proper moment map 
$$\mu:\  {\cal{M}}_H \rightarrow R.$$ In fact $\mu= ||\Phi||^2_{L^2}$, the $L^2$-norm of $\Phi$. The fixed point set of the action which is identical to the critical set of $\mu$ consists of  

(i). When $\Phi=0$, the moduli space ${\cal{M}}_u$ of flat $G$-bundles for some $u\in Z(G)$; 

(ii). The components of moduli spaces $\{ {\cal{M}} \}$ of representations of $\pi_1(S)$ in subgroups of $G$. This is when $\Phi\neq 0$. 

For $G=SU(2)$, then (2) just consists of the symmetric products of $S$. Let $z=\sqrt{-1}s$ be a complex number with $s>0$. Let $\omega_H$ be the natural symplectic form on ${\cal{M}}_H$ [Hi]. We then have the following formula of [Liu1]:

$$\int_{{\cal{M}}_H}e^{\omega_H+\sqrt{-1}z\mu}= \int_{{\cal{M}}_u}\frac{e^{\omega_u}}{e(z)}+\sum_{{\cal{M}}}\int_{{\cal{M}}}\frac{i_{{\cal{M}}}^* (e^{\omega_H+\sqrt{-1}z \mu})}{e_{\cal M}(z)}\ \ \ (4).$$
Here $i^*_{{\cal{M}}}$ denotes the restriction map, $e(z)$ and $e_{{\cal{M}}}(z)$ are respectively the equivariant Euler classes of the normal bundle of ${\cal{M}}_u$ and ${\cal{M}}$ in ${\cal{M}}_H$. Note that the terms on the right hand side are all computable by our previous results. See [Liu1] for detail and a more general formula.

This gives us an inductive formula to compute integrals on the Higgs moduli spaces, which may be useful for computations on the Teichmuller space of $S$ in connection with the work of [Go]. We hope to come back to this problem sometime later.

\vspace{.2in} 

{\bf 7. Invariants of knots and $3$-manifolds.} 

We first consider knots. Let $S^1\hookrightarrow S^3$ be a knot. The complement of the knot can be pushed to a $2$-complex $K$ with the homology of a circle. Consider a cell decomposition of $K$. Following [Mi], we can assume $K$ has only one $0$-cell.

Let $P$ be a flat principal $G$-bundle on $K$ with $G$ a semi-simple simply connected compact Lie group. Assume $K$ has $m$ two-cells, then $K$ has $m+1$ one-cells. The orientation on $K$ induces orientation on each cell. Let $V,\, E, \, F$ denote the space of $0$-cell, $1$-cells and $2$-cells respectively. Let 
$$A(P)= \mbox{Maps}(E, G)\simeq G^{m+1}, \ \ {\cal{G}}(P)= \mbox{Maps}(V, G)\simeq G$$
We consider elements in $A(P)$ as parallel transports along $1$-cells. For $x\in A(P)$, let $z\in {\cal{G}}(P)\simeq G$ acts on $x$ by 

$$x\rightarrow zx z^{-1}.$$

Note that this is precisely the ``lattice model'' of physicists as described in [Fo], pp41 and used in previous sections and [Liu]. Let $f$ be the holonomy map from $A(P)$ to $\mbox{Maps}(F, G)$. The boundary of $K$, $\partial K$ is a $1$-cell. We fix an element $c\in G$ and associate it to $\partial K$. In fact it is easy to see that we have flat connections on $P$ with fixed holonomy $c$ along the boundary of $K$. Then the holonomy map $f$ can be considered as a map from $G^m$ to $G^m$.

Consider the function on $\mbox{Maps}(F, G)\simeq G^m$ given by 
$$H_F(t, x, y)=\prod_j H(t, x_j, y_j)$$
where the product is over the $2$-cells $F$, and $x_j$, $y_j$ are respectively the components of $x,\ y\in \mbox{Maps}(F, G)\simeq G^m$.

Now consider the integral of the pull-back of $H_F$ by $f$,

$$I(t, c)= \int_{G^m}H_F(t, f(h), e)d\mbox{vol}\ \ \mbox{with} \ h\in \mbox{Maps}(E, G)\simeq G^m.$$where $d\mbox{vol}$ is the volume form on $G^m$ induced from the Killing form and $e\in G^m$ is the identity element.

When $t\rightarrow 0$, the same method as above gives us  

$$\mbox{lim}_{t\rightarrow 0}I(t, c)= \frac{1}{\#Z(G)}{\mbox{Vol}}(G)^m{\mbox{Vol}}(T)\int_{{\cal M}_c} d\nu$$
Where $${\cal M}_c=\{ f(h)=e \}/G\simeq {\mbox{Hom}}_c(\pi_1(K), G)/G $$is the moduli space of flat connecions on $P$ with fixed holonomy $c\in G$ along the boundary $\partial K$, and $d\nu$ which is a natural measure on ${\cal M}$ is the torsion of 

$$0\rightarrow {\cal{Z}}_c\stackrel{d\gamma}{\rightarrow} {\cal{G}}^m \stackrel{df}{\rightarrow} {\cal{G}}^m\rightarrow 0$$ 
which is precisely the cochain complex $(C^*(K, \partial K, \mbox{ad}\, P),\ \partial)$.

So we get 

$$\int_{{\cal M}_c}\tau( K, \partial K, \mbox{ad}\, P)=\frac{\#Z(G)}{{\mbox{Vol}}(T){\mbox{Vol}}(G)^m}\mbox{lim}_{t\rightarrow 0}\int_{G^m}H_F(t, f(h), e)d\mbox{vol}\ \ \ (5).$$ 
Usually ${\cal M}_c$ consists of isolated points and the left hand side is just the sum of Reidemeister torsions. By [Mi], the torsion is closely related to the Alexander knot polynomial. In equality (5), the right hand side usually is easy to compute and the left hand side can be considered as a generalization of knot invariants.

We can also consider the integral twisted by some class function $ A(h)$ on ${\mbox{Maps}}(E, G)\simeq G^m$,

$$\int_{G^m}H_F(t, f(h), c)A(h)d \mbox{vol}$$ which will give us different kind of invariant of the complex $K$, therefore the knot. The most interesting choice for $A(h)$ is the character of a representation of $G$ evaluated at the holonomy around some loop in $K$.

Now we turn to $3$-manifolds. We consider a closed $3$-manifold $M$. Let $P$ be a principal $G$-bundle on $M$. By constructing the ``lattice model'' on the $2$-skeleton of a cell decomposition of $M$, we are in the same situation as above. For simplicity we still assume there is only one $0$-cell. Let $A(h)$ be a class function on ${\mbox{ Maps}} (E, G)\simeq G^m$ where $m$ is the number of $1$-cells. Let $n$ be the number of $2$-cells, that is, $\mbox{Maps}(F, G)\simeq G^n$. Here let us only consider the simple situation. Assume the moduli space 

$$ {\cal M}={\mbox{Hom}}(\pi_1(M), G)/G $$ consists of finite number of isolated points which we denote by $\{ \lambda_j\}$. Let $\mbox{ad}\, P_j$ be the adjoint bundle associated to $P$ by $\lambda_j$, and let $\tau(M, \mbox{ad}\, P_j)$ denote the corresponding Reidemeister torsion. Let $f$ be the holonomy map 
$$f:\  \ \mbox{Maps}(E, G)\rightarrow \mbox{Maps}(F, G).$$
 For simplicity we omit the computational details, and just write down the final formula for this special case:

$$\sum_{\lambda_j\in {\cal M}}\tau(M, \mbox{ad}\, P_j)^{\frac{1}{2}} A(\lambda_j)=\frac{\# Z(G)}{{\mbox{Vol}}(G)^{n+1}}\mbox{lim}_{t\rightarrow 0}\int_{G^m}H_F(t, f(h), e)A(h)d\mbox{vol}\ \ \ $$
where the sum on the left hand side is over all equivalent classes of representations of $\pi_1(M)$ in $G$. The right hand side is easy to compute in examples. One interesting choice for $A(h)$ is the class function corresponding to the Chern-Simons functional. See Atiyah's lecture in [Ox].
 
These invariants for knots and $3$-manifolds are clearly related to the Jones-Witten type invariants, or the Casson type invariant [Jo], [Wa]. In fact our method is basically a finite dimensional version of path-integral as pointed out by Atiyah in [Ox]. Similar type invariants have been studied in [Jo] from simplicial point of view.

Note that, in general the left hand side of the above formula should be replaced by the integral over the moduli space ${\cal M}$ of flat connections on $P$. In this case the square root of torsion is a natural measure on ${\cal M}$. See [M]. As pointed by the referee, we have to deal with the singularities of the moduli spaces.

\vspace{.3in}

K. Liu, {\small Department of Mathematics, Stanford University, Stanford, CA 94305.}

\end{document}